\documentclass[fleqn]{mat01}
\usepackage{times,mathtimy,amssymb,latexsym}%,sidecap,graphicx,rangecite}
\begin{document}

\setcounter{page}{387} \firstpage{387}

\newtheorem{theore}{\bf Theorem}
\renewcommand\thetheore{\arabic{section}.\arabic{theore}}
\newtheorem{definit}[theore]{\rm DEFINITION}
\newtheorem{theor}[theore]{\bf Theorem}
\newtheorem{propo}[theore]{\rm PROPOSITION}
\newtheorem{lem}[theore]{\it Lemma}
\newtheorem{rem}[theore]{\it Remark}

\newtheorem{coro}[theore]{\rm COROLLARY}
\newtheorem{probl}[theore]{\it Problem}
\newtheorem{exampl}[theore]{\it Example}
\newtheorem{pot}[theore]{\it Proof of Theorem}

\renewcommand\theequation{\thesection\arabic{equation}}

\def\nota{\trivlist \item[\hskip \labelsep{\it Notations.}]}
\def\tog{\trivlist \item[\hskip \labelsep{\bf Theorem of Gadjiev.}]}

\def\d{\mbox{\rm d}}
\def\e{\mbox{\rm e}}

\title{Approximation of functions of two variables by certain linear
positive operators}

\markboth{Fatma Ta\c{s}delen et al}{Linear positive operators}

\author{FATMA TA\c{S}DELEN$^{*}$, ALI OLGUN$^{\dagger}$ and
G\"{U}LEN~BA\c{S}CANBAZ-TUNCA$^{*}$}

\address{$^{*}$Department of Mathematics, Faculty of Science, Ankara
University,
Tandogan~06100, Ankara, Turkey\\
\noindent$^{\dagger}$Department of Mathematics, Faculty of Science
and Arts, Kirikkale University, Yah\c{s}ihan~71450, Kirikkale, Turkey\\
\noindent E-mail: tasdelen@science.ankara.edu.tr;
aolgun@kku.edu.tr; tunca@science.ankara.edu.tr}

\volume{117}

\mon{August}

\parts{3}

\pubyear{2007}

\Date{MS received 24 February 2006}

\begin{abstract}
We introduce certain linear positive operators and study some
approximation properties of these operators in the space of
functions, continuous on a compact set, of two variables. We also
find the order of this approximation by using modulus of
continuity. Moreover we define an $r$th order generalization of
these operators and observe its approximation properties.
Furthermore, we study the convergence of the linear positive
operators in a weighted space of functions of two variables and
find the rate of this convergence using weighted modulus of
continuity.
\end{abstract}


\maketitle

\section{Introduction}

Let $f\in C([0,1])$. The well-known Bernstein polynomial of degree
$n$, denoted by $B_{n}(f;x)$ is
\begin{equation*}
B_{n}(f;x) :=\sum\limits_{k=0}^{n}p_{n,k}(x)
f\left(\frac{k}{n}\right),\quad n\in \mathbb{N=}\{1,2,\dots\},
\end{equation*}%
where
\begin{equation}
p_{n,k}(x) = \begin{pmatrix}{n}\\[.1pc] {k}\end{pmatrix}x^{k}(1-x)^{n-k}
\end{equation}%
and $x\in \lbrack 0,1]$ \cite{4}.

Let $x\in \lbrack 0,\infty )$ and $f\in C([0,\infty ))$.
Szazs--Mirakyan operators, denoted by $S_{n}(f;x)$ are
\begin{equation*}
S_{n}(f;x) :=\sum\limits_{k=0}^{\infty} q_{n,k}(x)
f\left(\frac{k}{n}\right),\quad n\in \mathbb{N},
\end{equation*}%
where
\begin{align}
q_{n,k}(x) = \e^{-nx}\frac{(nx) ^{k}}{k!}.
\end{align}
In \cite{3} approximation properties of $S_{n}(f;x)$ in weighted
spaces were studied. Some works by Szazs--Mirakyan or modified
Szazs--Mirakyan operators may be found in \cite{12,7,18} and
references therein.

Stancu \cite{16} introduced the following generalization of the
Bernstein polynomials. Let $f\in C([0,1])$. Stancu operators,
denoted by $(P_{n}^{(\alpha ,\beta ) }f)$ are
\begin{equation*}
(P_{n}^{(\alpha ,\beta ) }f) :=\sum\limits_{k=0}^{n}p_{n,k}(x) f\left(\frac{k+\alpha }{%
n+\beta }\right),\quad n\in \mathbb{N},
\end{equation*}
where $p_{n,k}(x)$ are the polynomials given by (1.1), $\alpha,
\beta$ are positive real numbers satisfying $0\leq \alpha \leq
\beta$.

Taking the operators, given above, into account we now introduce certain
linear positive operators of functions of two variables as follows:

Let $f\in C(\mathcal{R}), \mathcal{R}:=[0,1]\times \lbrack
0,\infty )$ and let the linear positive operators
$L_{m,n}^{\alpha_{i},\beta_{j}},\ j=1,2,\ $be defined as follows:
\begin{equation}
L_{m,n}^{\alpha_{i},\beta_{j}}:=\sum\limits_{k=0}^{\infty
}\sum\limits_{\nu =0}^{m}p_{m,\nu }(x) q_{n,k}(y) f\left(\frac{\nu
+\alpha_{1}}{m+\beta_{1}},\frac{k+\alpha_{2}}{n+\beta_{2}}\right)
\end{equation}
for $(x,y) \in \mathcal{R}, m,n\in \mathbb{N}, 0\leq \alpha
_{j}\leq \beta_{j},\ j=1,2,$ where $p_{m,\nu }(x)$ and
$q_{n,k}(y)$ are given in (1.1) and (1.2), respectively. In the
sequel, whenever we mention\ the operators
$L_{m,n}^{\alpha_{i},\beta_{j}}, j=1,2,$ it will be mentioned that
these are the operators given in (1.3). We use the notation
$\mathcal{R}_{A}$ to denote the following closed and bounded
region in $\mathbb{R}^{2}$,
\begin{align}
\mathcal{R}_{A}:=[0,1]\times \lbrack 0,A],\quad A>0.
\end{align}
In this paper we first study some approximation properties of the
sequence of linear positive operators given by (1.3) in the space
of functions, continuous on $\mathcal{R}_{A}$, and find the order
of this approximation using modulus of continuity. Moreover we
define an $r$th order generalization of $L_{m,n}^{\alpha
_{i},\beta_{j}}, j=1,2,$ on $\mathcal{R}_{A}$ extending the
results of Kirov \cite{14} and Kirov--Popova \cite{15} to the
linear positive operators\ $L_{m,n}^{\alpha_{i},\beta_{j}},
j=1,2,$ of functions of two variables and study its approximation
properties. The $r$th order generalization of some kind of linear
positive operators may also be found in \cite{1,9}.

We finally investigate the convergence of the sequence of linear
positive operators $L_{m,n}^{\alpha_{i},\beta_{j}}, j=1,2,$
defined on a weighted space of functions of two variables and find
the rate of this convergence by means of weighted modulus of
continuity.

If we take $p_{n,k}(y), k=0,1,\dots,n$ in place of $q_{n,k}(y)$ in
(1.3), then the operators $L_{m,n}^{\alpha_{i},\beta_{j}}, j=1,2,$
reduce to the generalized Bernstein polynomials of two variables
which were studied in \cite{5}.

Approximation of functions of one or two variables by some
positive linear operators in weighted spaces may be found in
\cite{8,9,13,17,18}.\newpage

\section{Preliminaries}

In this section we give some basic definitions which we shall use.
We denote by $\rho$ the function, continuous and satisfying $\rho
(x,y) \geqslant 1$ for $(x,y) \in \mathcal{R}$ and $\lim_{\vert
r\vert \rightarrow \infty }\rho (x,y) =\infty, r=(x,y) \cdot \rho$
is called a weight function. Let $B_{\rho }$ denote the set of
functions of two variables defined on $\mathcal{R}$ satisfying
$\vert f(x,y) \vert \leq M_{f}\rho (x,y),$ where $M_{f}>0$ is a
constant depending on $f$, and $C_{\rho }$ denote the set of
functions belonging to $B_{\rho },$ and continuous on
$\mathcal{R}.$ Clearly $C_{\rho }\subset B_{\rho }\cdot B_{\rho }$
and $C_{\rho }$ are called weighted spaces with norm $\Vert
f\Vert_{\rho }=\sup_{(x,y) \in \mathcal{R}}\frac{\vert f(x,y)
\vert }{\rho (x,y)}$, \cite{10,11}.

The Lipschitz class $\hbox{Lip}_{M}(\gamma)$ of the functions of
$f$ of two variables is given by \setcounter{equation}{0}
\begin{equation}
\vert f(x_{1},y_{1}) -f(x_{2},y_{2}) \vert \leq M[ (x_{1}-x_{2})
^{2}+(y_{1}-y_{2}) ^{2}] ^{\frac{\gamma }{2}},
\end{equation}
$(x_{1},y_{1}), (x_{2},y_{2}) \in \mathcal{R},$ where $M>0,
0<\gamma \leq 1$ and $f\in C(\mathcal{R})$. The full modulus of
continuity of $f \in C(\mathcal{R}_{A}),$ denoted by
$w(f;\delta),$ is defined as follows:
\begin{equation}
w(f;\delta ) =\max\limits_{\sqrt{(x_{1}-x_{2}) ^{2}+(y_{1}-y_{2})
^{2}}\leq \delta }\vert f(x_{1},y_{1}) -f(x_{2},y_{2}) \vert.
\end{equation}
Partial modulus of continuity with respect to $x$ and $y$ are
given by
\begin{equation}
w^{(1) }(f;\delta ) =\max\limits_{0\leq y\leq A}\ \
\max\limits_{\vert x_{1}-x_{2}\vert \leq \delta }\vert f(x_{1},y)
-f(x_{2},y) \vert
\end{equation}
and
\begin{equation}
w^{(2) }(f;\delta ) =\max\limits_{0\leq x\leq 1}\ \
\max\limits_{\vert y_{1}-y_{2}\vert \leq \delta }\vert f(x,y_{1})
-f(x,y_{2}) \vert,
\end{equation}
respectively. We shall also need the following properties of the
full and partial modulus of continuity
\begin{equation}
w(f;\lambda \delta ) \leq (1+[ \lambda ] ) w(f;\delta )
\end{equation}
for any $\lambda$. Here $[\lambda]$ is the greatest integer that
does not exceed $\lambda$. Moreover, it is known that when $f$ is
uniformly continuous, then $\lim_{\delta \rightarrow 0}w(f;\delta
) =0$ and
\begin{equation}
\vert f(t,\tau ) -f(x,y) \vert \leq w(f;\sqrt{(t-x) ^{2}+(\tau -y)
^{2}}),
\end{equation}
$(t,\tau ), (x,y) \in \mathcal{R}_{A}.$ The analogous properties
are satisfied by the partial modulus of continuity.

\section{Lemmas and theorems on \pmb{$\mathcal{R}_{A}$}}

In this section we give some classical approximation properties of
the operators $L_{m,n}^{\alpha_{i},\beta_{j}}, j=1,2,$ on the
compact set $\mathcal{R}_{A}$.

\begin{lem}
Let $\alpha_{j},\beta_{j}, j=1,2,$ be the fixed positive numbers
such that $0\leq \alpha_{j}\leq \beta_{j}$. Then we have
\begin{align*}
L_{m,n}^{\alpha_{i},\beta_{j}}(1;x,y) =1,%\\[.4pc]
\end{align*}
\begin{align*}
L_{m,n}^{\alpha_{i},\beta_{j}}(t;x,y) &=\frac{mx+\alpha_{1}}{m+\beta_{1}},\\[.4pc]
L_{m,n}^{\alpha_{i},\beta_{j}}(\tau ;x,y) &=\frac{nx+\alpha
_{2}}{n+\beta_{2}}, \\[.4pc]
L_{m,n}^{\alpha_{i},\beta_{j}}(t^{2}+\tau ^{2};x,y) &=\frac{%
(m^{2}-m) x^{2}+(2\alpha_{1}+1) mx+\alpha_{1}^{2}}{%
(m+\beta_{1}) ^{2}} \\[.4pc]
&\quad\, +\frac{n^{2}y^{2}+(2\alpha_{2}+1) ny+\alpha
_{2}^{2}}{(n+\beta_{2}) ^{2}}
\end{align*}%
for all $m,n\in N.$
\end{lem}

Taking (3.1) into account we now give the following Baskakov type
theorem (see \cite{2} to get the approximation to $f(x,y) \in
C(\mathcal{R}_{A}),$ satisfying $\vert f(x,y) \vert \leq
M_{f}(1+x^{2}+y^{2}),$ by $L_{m,n}^{\alpha_{i},\beta_{j}}, j=1,2.$

\begin{theore}[\!]
Let $f(x,y) \in C(\mathcal{R}_{A})$ and $\vert f(x,y) \vert \leq
M_{f}(1+x^{2}+y^{2})$ for $(x,y) \in R$. Here $M_{f}$ is a
constant depending on $f$. Then $\Vert
L_{m,n}^{\alpha_{i},\beta_{j}}(f;x,y) -f(x,y) \Vert
_{C(\mathcal{R}_{A}) }\rightarrow 0${\rm ,}  as  $m,n \rightarrow
\infty$ if and only if \setcounter{equation}{0}
\begin{align}
\Vert L_{m,n}^{\alpha_{i},\beta_{j}}(1;x,y) -1\Vert
_{C(\mathcal{R}_{A}) } &\rightarrow 0,\nonumber\\[.4pc]
\Vert L_{m,n}^{\alpha_{i},\beta_{j}}(t;x,y) -x\Vert
_{C(\mathcal{R}_{A}) } &\rightarrow 0,\nonumber\\[.4pc]
\Vert L_{m,n}^{\alpha_{i},\beta_{j}}(\tau ;x,y) -y\Vert
_{C(\mathcal{R}_{A}) } &\rightarrow 0,\nonumber\\[.4pc]
\Vert L_{m,n}^{\alpha_{i},\beta_{j}}(t^{2}+\tau ^{2};x,y)
-(x^{2}+y^{2}) \Vert_{C(\mathcal{R}_{A}) } &\rightarrow 0,
\end{align}
as $m,n\rightarrow \infty$ for $(x,y) \in R_{A}.$
\end{theore}

\begin{proof}
Since the necessity is clear, then we need only to prove the
sufficiency. Let \hbox{$(t,\tau ), (x,y) \in \mathcal{R}_{A}.$} By
the uniform continuity of $f$ on $\mathcal{R}_{A}$ we get that
for~each~\hbox{$\varepsilon >0$} there exists a number $\delta >0$
such that $\vert f(t,\tau ) -f(x,y) \vert <\varepsilon,$ whenever
$\sqrt{(t-x) ^{2}+(\tau -y) ^{2}}<\delta$. Now let $(x,y) \in
\mathcal{R}_{A}$ and $(t,\tau ) \in \mathcal{R}$ and let
$(x_{1},y_{1})$ be an arbitrary boundary point of
$\mathcal{R}_{A}$ such that $0\leq x_{1}\leq 1, 0\leq y_{1}\leq
A$. Since $f$ is continuous on the boundary points also, then for
each $\varepsilon >0$ there exists a $\delta >0$ such that
\begin{equation*}
\vert f(t,\tau ) -f(x,y) \vert \leq \vert f(t,\tau )
-f(x_{1},y_{1}) \vert +\vert f(x_{1},y_{1}) -f(x,y) \vert
<\varepsilon
\end{equation*}
whenever $\sqrt{(t-x) ^{2}+(\tau -y) ^{2}} <\delta$. Finally let
$(x,y) \in \mathcal{R}_{A}$ and $(t,\tau ) \in \mathcal{R}$ and
let $\sqrt{(t-x) ^{2}+(\tau -y) ^{2}}>\delta$. Then easy
calculations show that
\begin{align*}
\vert f(t,\tau ) -f(x,y) \vert &\leq M_{f}((t-x) ^{2}+(\tau -y)
^{2}) \left(\frac{2}{\delta ^{2}}+2+\frac{3}{\delta
^{2}}(x^{2}+y^{2})\right)\\[.4pc]
&\leq C\left(\frac{(t-x) ^{2}+(\tau -y) ^{2}}{\delta ^{2}}\right).
\end{align*}
Here\ $C>0$ is a constant. Therefore we get
\begin{equation}
\vert f(t,\tau ) -f(x,y) \vert \leq
\varepsilon + C\left(\frac{(t-x) ^{2}+(\tau -y) ^{2}%
}{\delta ^{2}}\right),
\end{equation}
for $(t,\tau ) \in \mathcal{R}, (x,y) \in \mathcal{R}_{A}.$
Applying $L_{m,n}^{\alpha_{i},\beta_{j}}$ to (3.2) we get
\begin{align*}
\vert L_{m,n}^{\alpha_{i},\beta_{j}}(f(t,\tau ) ;x,y) -f(x,y)
\vert &\leq L_{m,n}^{\alpha_{i},\beta_{j}}(\vert f(t,\tau )
-f(x,y)\vert ;x,y) \\[.4pc]
&\quad\,+\Vert f\Vert \vert L_{m,n}^{\alpha_{i},\beta_{j}}(1;x,y)
-1\vert.
\end{align*}
Using (3.2) in the last inequality and taking (3.1) into account,
sufficiency is obtained easily. \hfill $\Box$
\end{proof}

We note that if we take $f(x,y)$ to be bounded on $\mathbb{R}^{2}$
in the previous theorem, then we easily obtain that $\Vert
L_{m,n}^{\alpha_{i},\beta_{j}}(f;x,y) -f(x,y) \Vert
_{C(\mathcal{R}_{A}) }\rightarrow 0,$ as $m, n\rightarrow \infty$
satisfied from Lemma~3.1 by analogous Korovkin's theorem proved by
Volkov \cite{19}.

The following theorem gives the rate of convergence of the
sequence of linear positive operators $\{L_{m,n}^{\alpha
_{i},\beta_{j}}\}$ to $f,$ by means of partial and full modulus of
continuity.

\begin{theore}[\!]
Let $f\in C(\mathcal{R}_{A})$. Then the following inequalities
\begin{align}
\hskip -4pc {\rm (a)} \hskip 2.98pc &\Vert
L_{m,n}^{\alpha_{i},\beta _{j}}(f;x,y) -f(x,y)
\Vert_{C(\mathcal{R}_{A}) }\leq \frac{3}{2}
\{ w^{(1) }(f;\delta_{m}) +w^{(2) }(f;\delta_{n}) \},\\[.4pc]
\hskip -4pc {\rm (b)} \hskip 2.98pc  &\Vert L_{m,n}^{\alpha
_{i},\beta_{j}}(f;x,y) -f(x,y) \Vert_{C(\mathcal{R}_{A}) }\leq
\frac{3}{2}w(f;\delta_{m,n})
\end{align}
hold{\rm ,} where $R_{A}$ is the closed and bounded region given
by $(1.4)$. $w^{(1)}, w^{(2)}$ and $w$ are given by $(2.3), (2.4)$
and $(2.2)$ respectively{\rm ,} and $\delta_{m}, \delta_{n},
\delta_{m,n}$ are
\begin{equation}
\delta_{m}=\frac{\sqrt{4\beta_{1}^{2}+m}}{m+\beta_{1}},\quad
\delta _{n}=\frac{\sqrt{\beta_{2}^{2}A^{2}+nA}}{n+\beta_{2}},\quad
\delta _{m,n}=\sqrt{\delta_{m}^{2}+4\delta_{n}^{2}},
\end{equation}%
respectively.
\end{theore}

\begin{proof}
From (1.3) we have
\begin{align}
\vert L_{m,n}^{\alpha_{i},\beta_{j}}(f(t,\tau ) ;x,y) -f(x,y)
\vert &\leq \e^{-ny}\sum\limits_{k=0}^{\infty }\sum\limits_{\nu
=0}^{m}\begin{pmatrix}{m}\\[.1pc]{\nu }\end{pmatrix}x^{\nu }(1-x) ^{m-\nu }\frac{(ny)
^{k}}{k!}\nonumber\\[.45pc]
&\quad\times \left\vert f\left(\frac{\nu +\alpha_{1}}{m+\beta_{1}},\frac{%
k+\alpha_{2}}{n+\beta_{2}}\right) -f(x,y) \right\vert . %\TCItag{3.6}
\end{align}
Let us first add and drop the function $f\big(\frac{\nu +\alpha
_{1}}{m+\beta_{1}},y\big)$ inside the absolute value sign on the
right-hand side of (3.6). Using the analogous property of (2.6)
for the partial modulus of continuity and finally applying the
Cauchy--Schwartz inequality to the resulting term, then we arrive
at (3.3) on $\mathcal{R}_{A},$ which proves (a). Using (2.6)
directly in (3.6) and applying Cauchy--Schwartz inequality to the
resulting term we then reach to (3.4) on $\mathcal{R}_{A},$ which
gives (b).\hfill $\Box$
\end{proof}

\begin{coro}\label{abelian}$\left.\right.$\vspace{.5pc}

\noindent Let $f\in {\rm Lip}_{M}(\gamma)$. Then the inequality
\begin{equation*}
\vert L_{m,n}^{\alpha_{i},\beta_{j}}(f;x,y) -f(x,y) \vert \leq
M_{1}'\delta_{m,n}^{\gamma }
\end{equation*}%
holds{\rm ,} where $M_{1}'=\frac{3}{2}M_{1}, M_{1}>0$ and $\delta
_{m,n}$ is given in $(3.5)$.
\end{coro}

\begin{coro}\label{abelian}$\left.\right.$\vspace{.5pc}

\noindent  If $f$ satisfies the following Lipschitz conditions
\begin{equation*}
\vert f(x_{1},y) -f(x_{2},y) \vert \leq M_{2}\vert
x_{1}-x_{2}\vert ^{\alpha }
\end{equation*}
and
\begin{equation*}
\vert f(x,y_{1}) -f(x,y_{2}) \vert \leq M_{3}\vert
y_{1}-y_{2}\vert ^{\beta },
\end{equation*}
$0<\alpha, \beta \leq 1, M_{j}>0, j=2,3,$ then the inequality
\begin{equation*}
\vert L_{m,n}^{\alpha_{i},\beta_{j}}(f;x,y) -f(x,y) \vert \leq
M_{2}'\delta_{m}^{\alpha }+M_{3}'(2\delta_{n}) ^{\beta }
\end{equation*}
holds{\rm ,} where $M_{2}' =\frac{3}{2}M_{2}$ and $M_{3}' =
\frac{3}{2}M_{3},$ and $\delta_{m}, \delta_{n}$ are given in
$(3.5)$.
\end{coro}

\section{A generalization of order \pmb{$r$} of
\pmb{$L_{m,n}^{\alpha_{i},\beta_{j}}$}}

Let $C^{r}(\mathcal{R}_{A}) $, $r\in \mathbb{N}\cup \{ 0\},$
denote the set of all functions\ $f$ having all continuous partial
derivatives up to order $r$ at $(x,y) \in \mathcal{R}_{A}.$

By $(L_{m,n}^{\alpha_{i},\beta_{j}})^{[ r]},j=1,2,$ we denote the
following generalization of $L_{m,n}^{\alpha_{i},\beta_{j}}$:
\setcounter{equation}{0}
\begin{align}
(L_{m,n}^{\alpha_{i},\beta_{j}}) ^{[ r] }(f;x,y) &:=
\e^{-ny}\sum\limits_{k=0}^{\infty }\sum\limits_{\nu =0}^{m}
\begin{pmatrix}{m}\\[.1pc]{\nu }\end{pmatrix}x^{\nu }(1-x) ^{m-\nu }\frac{(ny) ^{k}}{k!}\nonumber\\[.4pc]
&\quad\, \times P_{r,\big(\frac{\nu +\alpha_{{1}}}{m+\beta
_{1}},\frac{k+\alpha_{{2}}}{n+\beta_{2}}\big) }\left(x-\frac{\nu
+\alpha_{{1}}}{m+\beta_{1}},y-\frac{k+\alpha_{{2}}}{n+\beta
_{2}}\right),\hskip -1pc\phantom{0}
\end{align}
$0\leq \alpha_{j}\leq \beta_{j}$, $j=1,2,$ where
\begin{align}
&P_{r,\big(\frac{\nu +\alpha_{{1}}}{m+\beta_{1}},\frac{k+\alpha_{{2}}}{%
n+\beta_{2}}\big) }\left(x-\frac{\nu +\alpha_{{1}}}{m+\beta_{1}},y-%
\frac{k+\alpha_{{2}}}{n+\beta_{2}}\right)\nonumber\\[.6pc]
&\quad= \sum\limits_{h=0}^{r}\sum
\limits_{i+j=h}\frac{1}{h!}\begin{pmatrix}{h}\\[.1pc]{j}\end{pmatrix}f_{x^{i}y^{j}}\left(\frac{\nu
+\alpha_{{1}}}{m+\beta_{1}},\frac{k+\alpha_{{2}}}{n+\beta
_{2}}\right)\nonumber\\[.6pc]
&\qquad\ \times \left[ x-\frac{\nu +\alpha_{{1}}}{m+\beta
_{1}}\right] ^{i}\left[ y-\frac{k+\alpha_{{2}}}{n+\beta
_{2}}\right]^{j},
\end{align}
$f_{x^{i}y^{j}}$ denotes the partial derivatives of $f$, i.e.:
$f_{{x^{i}y^{j}}}:=\frac{\partial^{r}}{\partial x^{i}\partial
y^{j}}f(x,y)$.

$(L_{m,n}^{\alpha_{i},\beta_{j}}) ^{[ r] }$ are called the $r$th
order of $L_{m,n}^{\alpha_{i},\beta_{j}}$ (see \cite{14} for one
variable). Obviously $(L_{m,n}^{\alpha_{i},\beta_{j}}) ^{[ r] }$
reduce to $L_{m,n}^{\alpha_{i},\beta_{j}},$ when $r=0$.

Now let us write
\begin{equation}
\left(x-\frac{\nu +\alpha_{{1}}}{m+\beta_{1}},y-\frac{k+\alpha_{{2}}}{%
n+\beta_{2}}\right) =u(\alpha ,\beta ) ,
\end{equation}
where $(\alpha ,\beta ) $ is a unit vector, $u>0$ and let
\begin{align}
F(u)  &=f\left(\frac{\nu +\alpha_{{1}}}{m+\beta_{1}} +u\alpha
,\frac{k+\alpha_{{2}}}{n+\beta_{2}}+u\beta\right)\nonumber\\[.65pc]
&=f\left[ \frac{\nu +\alpha_{{1}}}{m+\beta_{1}}+\left(x-\frac{\nu
+\alpha_{{1}}}{m+\beta_{1}}\right) ,\frac{k+\alpha_{{2}}}{n+\beta_{2}}%
+\left(y-\frac{k+\alpha_{{2}}}{n+\beta_{2}}\right)\right].
\end{align}
It is clear that Taylor's formula for $F(u)$ at $u=0$ turns into
Taylor's formula for $f(x,y) $ at $\big(\frac{\nu
+\alpha_{{1}}}{m+\beta_{1}},\frac{k+\alpha_{{2}}}{n+\beta
_{2}}\big)$. Morever $r$th derivative takes the form
\begin{equation}
F^{(r) }(u)
=\sum\limits_{i+j=r}\begin{pmatrix}{r}\\[.1pc]{j}\end{pmatrix}
f_{x^{i}y^{j}}\left(\frac{\nu +\alpha_{{1}}}{m+\beta_{1}}+u\alpha
,\frac{k+\alpha_{{2}}}{n+\beta_{2}}+u\beta\right) \alpha ^{i}\beta
^{j},
\end{equation}
$r\in \mathbb{N}$ (see Chapter~3 of \cite{6}).

By means of the modification stated above ((4.3)--(4.5)), we get
the following result.

\setcounter{theore}{0}
\begin{theore}[\!]
Let $f\in C^{r}(\mathcal{R}_{A})$ and $F^{(r) }(u) \in {\rm
Lip}_{M}(\gamma)$. Then the inequality
\begin{align}
&\Vert (L_{m,n}^{\alpha_{i},\beta_{j}}) ^{[ r] }(f;x,y)
-f(x,y)\Vert_{C(\mathcal{R}_{A}) }\nonumber\\[.6pc]
&\quad\leq \frac{\gamma M}{\gamma +r}\frac{B(\gamma ,r) }{(r-1)
!}\Vert L_{m,n}^{\alpha_{i},\beta_{j}}(\vert (x,y)-(t,\tau )\vert
^{r+\gamma };x,y) \Vert_{C(\mathcal{R}_{A})}
\end{align}
holds{\rm ,} where $F^{(r)}(u)$ are given by $(4.5)${\rm ,}
$B(\gamma ,r)$ is the well-known beta function{\rm ,} $r,m,n\in
N,\ 0<\gamma \leq 1$ and $M>0$.
\end{theore}

\begin{proof}
From (4.1) and (4.2) we have
\begin{align}
&f(x,y)-(L_{m,n}^{\alpha_{i},\beta_{j}}) ^{[ r] }(f;x,y)\nonumber\\[.45pc]
&\quad =\sum\limits_{\nu =0}^{m}\begin{pmatrix}{m}\\[.1pc]{\nu }\end{pmatrix}
x^{\nu }(1-x) ^{m-\nu }\e^{-ny\frac{(ny) ^{k}}{k!}%
}\nonumber\\[.45pc]
&\qquad\ \times \sum\limits_{k=0}^{\infty }
\left\{ f(x,y)-\sum\limits_{h=0}^{r}\frac{1}{h!}\sum\limits_{i+j=h}\begin{pmatrix}{h}\\[.1pc]{j}\end{pmatrix}
f_{x^{i}y^{j}}\right.\nonumber\\[.45pc]
&\qquad\ \times \left.\left(\frac{\nu +\alpha_{{1}}}{m+\beta
_{1}},\frac{k+\alpha_{{2}}}{n+\beta_{2}}\right) \left[ x-\frac{\nu
+\alpha_{{1}}}{m+\beta_{1}}\right] ^{i}\left[ y-\frac{k+\alpha
_{{2}}}{n+\beta_{2}}\right] ^{j}\right\}.
\end{align}

We now consider Taylor's formula with the remainder for the
functions of two variables. Using the integral form of the
remainder term that appeared in (4.7), we arrive at
\begin{align}
&f(x,y)-P_{r,\big(\frac{\nu +\alpha_{{1}}}{m+\beta
_{1}},\frac{k+\alpha_{{2}}}{n+\beta_{2}}\big) }\left(x-\frac{\nu
+\alpha_{{1}}}{m+\beta
_{1}},y-\frac{k+\alpha_{{2}}}{n+\beta_{2}}\right)\nonumber\\[.55pc]
&\quad =\frac{1}{(r-1)!}\int_{0}^{1}\sum\limits_{i+j=h}
\begin{pmatrix}{h}\\[.1pc]{j}\end{pmatrix}
\left[ x-\frac{\nu +\alpha_{{1}}}{m+\beta_{1}}\right] ^{i}\left[y-
\frac{k+\alpha_{{2}}}{n+\beta_{2}}\right] ^{j}\nonumber\\[.55pc]
&\qquad\ \times f_{x^{i}y^{j}}\left[ \frac{\nu +\alpha
_{{1}}}{m+\beta_{1}}+t\left(x-\frac{\nu +\alpha_{{1}}}{m+\beta
_{1}}\right) ,\frac{k+\alpha_{{2}}}{n+\beta
_{2}}+t\left(y-\frac{k+\alpha_{{2}}}{n+\beta_{2}}\right)
\right]\nonumber\\[.55pc]
&\qquad\ \times (1-t)^{r-1}\d t.
\end{align}
Taking (4.3)--(4.5) into account, (4.8) turns into the following
form:
\begin{equation}
F(u) -\sum\limits_{h=0}^{r}\frac{1}{h!}F^{(h) }(0)
u^{h}=\frac{u^{r}}{(r-1) !}\int_{0}^{1} [ F^{(r) }(tu) -F^{(r)
}(0) ] (1-t) ^{r-1}\d t.
\end{equation}
From (4.3), (4.8), (4.9) and the fact that $F^{(r) }\in
{\rm Lip}_{M}(\gamma ) $ it follows that%
\begin{align}
&f(x,y)-P_{r,\big(\frac{\nu +\alpha_{{1}}}{m+\beta
_{1}},\frac{k+\alpha_{{2}}}{n+\beta_{2}}\big) }\left(x-\frac{\nu
+\alpha_{{1}}}{m+\beta_{1}},y-\frac{k+\alpha_{{2}}}{n+\beta
_{2}}\right)\nonumber\\[.65pc]
&\quad=\left\vert F(u) -\sum\limits_{h=0}^{r}\frac{1}{h!}F^{(h) }(0) u^{h}\right\vert\nonumber\\[.6pc]
&\quad\leq \frac{\vert u\vert ^{r}}{(r-1) !}
\int_{0}^{1}[ F^{(r) }(tu) -F^{(r) }(0) ] (1-t) ^{r-1}\d t\nonumber\\[.6pc]
&\quad\leq \frac{\vert u\vert ^{r+\gamma }}{(r-1) !}MB(\gamma +1,r)\nonumber\\[.6pc]
&\quad\leq \frac{M}{(r-1) !}\frac{\gamma }{\gamma +r}B(\gamma, r) \vert u\vert ^{r+\gamma }\nonumber\\[.6pc]
&\quad\leq \frac{M}{(r-1) !}\frac{\gamma }{\gamma +r}B(\gamma ,r)
\left\vert x-\frac{\nu +\alpha_{{1}}}{m+\beta
_{1}},y-\frac{k+\alpha_{{2}}}{n+\beta_{2}}\right\vert ^{r+\gamma
}.
\end{align}
Hence combining (4.7) and (4.10), we obtain (4.6), which completes
the proof. \hfill $\Box$
\end{proof}

Now we take a function $g\in C(\mathcal{R}_{A})$ which is given by
\begin{equation}
g(t,\tau ) =\vert (x,y) -(t,\tau ) \vert ^{r+\gamma }.
\end{equation}
Obviously $g(x,y) =0.$ From Theorem~3.2 it follows that
\begin{equation*}
\Vert L_{m,n}(g;x,y)\Vert_{C(\mathcal{R}_{A}) }\rightarrow
0\quad\hbox{as}\quad m,n\rightarrow \infty .
\end{equation*}%
From (4.6) we arrive at the following result:
\begin{equation*}
\Vert (L_{m,n}^{\alpha_{i},\beta_{j}}) ^{[ r] }(f;x,y)
-f(x,y)\Vert_{C(\mathcal{R}_{A}) }\rightarrow 0\quad\hbox{as}\quad
m,n\rightarrow \infty .
\end{equation*}

Using (3.3) and Corollary~3.4 we get to the following results by
means of Theorem~4.1.

\begin{coro}\label{abelian}$\left.\right.$\vspace{.5pc}

\noindent Let $f\in C^{r}(\mathcal{R}_{A})$ and $F^{(r) }\in {\rm
Lip}_{M}(\gamma)$. Then the inequality
\begin{align*}
\Vert (L_{m,n}^{\alpha_{i},\beta_{j}}) ^{[ r] }(f;x,y)
-f(x,y)\Vert_{C(\mathcal{R}_{A})
}\leq \frac{MB(\gamma ,r) }{(r-1) !}\frac{\gamma }{%
\gamma +r}\frac{3}{2}w(g;\delta_{m,n})
\end{align*}
holds{\rm ,} where $F^{(r)}, \delta_{m,n}$ and $g$ are given by
$(4.5), (3.5)$ and $(4.11)${\rm ,} respectively.
\end{coro}

\begin{coro}\label{abelian}$\left.\right.$\vspace{.5pc}

\noindent Let $f\in C^{r}(\mathcal{R}_{A})$ and $F^{(r) }\in{\rm
Lip}_{M}(\gamma)${\rm ,} and assume that $g\in {\rm
Lip}_{(1+A^{2})^{\frac{r}{2}}}(\gamma)$ in Corollary~$3.4$. Then
we arrive at
\begin{align*}
\Vert (L_{m,n}^{\alpha_{i},\beta_{j}}) ^{[ r] }(f;x,y)
-f(x,y)\Vert_{C(\mathcal{R}_{A})
}\leq \frac{M(1+A^{2}) ^{\frac{r}{2}}}{(r-1) !}\frac{%
\gamma }{\gamma +r}B(\gamma ,r) \delta_{m,n}^{\gamma },
\end{align*}
where $\delta_{m,n}$ is given by $(3.5)$.
\end{coro}

\section{Weighted approximation of functions of two variables by
\pmb{$L_{m,n}^{\alpha_{i},\beta_{j}}$}}

In this section we investigate the convergence of the sequence
$\{L_{m,n}^{\alpha_{i},\beta_{j}}\}$ mapping the weighted space
$C_{\rho}$ into $B_{\rho_{1}}$. We also study the rates of
convergence of the sequence $\{ L_{m,n}^{\alpha_{i},\beta_{j}}\} $
defined on weighted spaces. In the rest of the article $\rho$ will
be given by $\rho (x,y) =$\break $1+x^{2}+y^{2}.$

We first give the following important Korovkin type theorem (in
weighted spaces) proved by Gadjiev in \cite{11}.

\begin{tog}
{\it Let $\{ A_{n}\}$ be the sequence of linear positive operators
mapping from $C_{\rho }(\mathbb{R}^{m})$ into $B_{\rho
}(\mathbb{R}^{m}), m\geqslant 1,$ and satisfying the conditions
\begin{align*}
\Vert A_{n}(1;\mathbf{x}) -1\Vert_{\rho } &\rightarrow 0,\quad
\Vert A_{n}(t_{j};\mathbf{x})-x_{j}\Vert_{\rho }\rightarrow 0,\ j=1,\dots,m,\\[.3pc]
\Vert A_{n}(\vert t\vert ^{2};\mathbf{x}) -\vert \mathbf{x}\vert
^{2}\Vert_{\rho } &\rightarrow 0
\end{align*}
as $n\rightarrow \infty$ for $\rho (\mathbf{x}) =1+\vert
\mathbf{x}\vert^{2},\ x\in R^{m}.$ Then there exists a function
$f^{\ast }\in C_{\rho }(\mathbb{R}^{m})$ such that $\Vert
A_{n}(f^{\ast };\mathbf{x}) -\mathbf{f}^{\ast }(\mathbf{x}) \Vert
_{\rho }\geqslant 1$.}\vspace{.5pc}
\end{tog}

By taking the result of the last theorem into account we conclude
that verifying the conditions of the above theorem by the
operators $L_{m,n}^{\alpha_{i},\beta_{j}}, j=1,2,$ is not
sufficient for $L_{m,n}^{\alpha_{i},\beta_{j}}$ to be convergent
to any function $f$ in $\rho$ norm. Hence we need to show the
convergence in another norm for any function in
$C_{\rho}(\mathcal{R})$. For this purpose, we now give the
following lemma, which we shall use.\pagebreak

\setcounter{theore}{0}
\begin{lem}
The operators $L_{m,n}^{\alpha_{i}, \beta_{j}}$ possess the
following{\rm :}
\begin{enumerate}
\renewcommand\labelenumi{{\rm (\alph{enumi})}}
\leftskip .15pc

\item $\{L_{m,n}^{\alpha_{i},\beta_{j}}\}, m,n\in N${\rm ,}
is the sequence of linear positive operators from the weighted
space $C_{\rho }(\mathcal{R})$ into the weighted space
$B_{\rho}(\mathcal{R})$.

\item The norms $\Vert L_{m,n}^{\alpha_{i},\beta_{j}}
\Vert_{C_{\rho }\rightarrow B_{\rho }}$ are uniformly bounded {\rm
(}i.e. there exists an $M>0$ such that $\Vert
L_{m,n}^{\alpha_{i},\beta_{j}}\Vert_{C_{\rho }\rightarrow B_{\rho
}}\leq M${\rm )}.
\end{enumerate}
\end{lem}

\begin{proof}
From Lemma~3.1 we easily obtain that
\begin{equation*}
\vert L_{m,n}^{\alpha_{i},\beta_{j}}(\rho ;x,y) \vert \leq
M(1+x^{2}+y^{2})
\end{equation*}
which proves (a). Taking Lemma~3.1 into account we get the
following inequality: \setcounter{equation}{0}
\begin{align}
&\Vert L_{m,n}^{\alpha_{i},\beta_{j}}\Vert_{C_{\rho }\rightarrow
B_{\rho }}\nonumber\\[.55pc]
&\quad\leq \Vert L_{m,n}^{\alpha_{i},\beta_{j}}(\rho
;x,y) -\rho (x,y) \Vert_{\rho }+1\nonumber\\[.55pc]
\hskip -4pc &\quad\leq \Vert L_{m,n}^{\alpha_{i},\beta_{j}}(1;x,y)
-1\Vert_{\rho
}+\Vert L_{m,n}^{\alpha_{i},\beta_{j}}(t^{2}+\tau ^{2};x,y) -x^{2}-y^{2}\Vert_{\rho }+1\nonumber\\[.45pc]
&\quad=\underset{(x,y) \in \mathcal{R}}{\sup }\left\vert
\frac{(m^{2}-m) x^{2}+(2\alpha_{1}+1)
mx+\alpha_{1}^{2}}{(m+\beta_{1})
^{2}}\right.\nonumber\\[.55pc]
&\left.\qquad\ +\frac{n^{2}y^{2}+(2\alpha_{2}+1) ny+\alpha
_{2}^{2}}{(n+\beta_{2}) ^{2}}-x^{2}-y^{2}\phantom{\underset{(x,y)
\in \mathcal{R}}{\sup }}\hskip -2.4pc\right\vert
\frac{1}{1+x^{2}+y^{2}}+1
\end{align}
so (b) is obtained from (5.1), which completes the proof.\hfill
$\Box$
\end{proof}

Now the following theorem shows the convergence of the sequence of
linear positive operators $\{ L_{m,n}^{\alpha_{i},\beta_{j}}\},$
mapping from $C_{\rho }$ into $B_{\rho_{1}},$ in $\rho_{1}$ norm.

\begin{theore}[\!]
Let $\rho_{1}(x,y)$ be a weight function satisfying
\begin{equation}
\underset{\vert x\vert \rightarrow \infty }{\lim }\frac{\rho (x,y)
}{\rho_{1}(x,y) }=0.
\end{equation}
Then $\Vert L_{m,n}^{\alpha_{i},\beta_{j}}(f;x,y)
-f(x,y)\Vert_{\rho_{1}}\rightarrow 0,\ m,n\rightarrow \infty$ for
all $f\in C_{\rho }(\mathcal{R})${\rm ,} where $x=(x,y) \in R$.
\end{theore}

\begin{proof}
Let us denote the region $[ 0,1] \times [ 0,s], s>0,$ by
$\mathcal{R}_{s}.$ Therefore we have
\begin{align}
&\Vert L_{m,n}^{\alpha_{i},\beta_{j}}(f;x,y)
-f(x,y)\Vert_{\rho_{1}}\nonumber\\[.55pc]
&\quad=\underset{(x,y) \in \mathcal{R}}{\sup }\left\vert
\frac{L_{m,n}^{\alpha_{i},\beta_{j}}(f;x,y)
-f(x,y)}{\rho_{1}(x,y) }\right\vert \nonumber\\[.55pc]
&\quad=\underset{(x,y) \in \mathcal{R}_{s}}{\sup }\frac{\vert
L_{m,n}^{\alpha_{i},\beta_{j}}(f;x,y) -f(x,y)\vert }{ \rho (x,y)
}\frac{\rho (x,y) }{\rho_{1}(x,y) }\nonumber\\[.55pc]
&\qquad\ +\underset{(x,y) \in \mathcal{R\smallsetminus R}_{s}}
{\sup }\frac{\vert L_{m,n}^{\alpha_{i},\beta_{j}}(f;x,y)
-f(x,y)\vert }{\rho (x,y) }\frac{\rho (x,y) }{ \rho_{1}(x,y) }.
\end{align}
Since ${\rho }/{\rho_{1}}$ is bounded on $\mathcal{R}_{s},$ the
first term on the right-hand side of (5.3) approaches zero when
$m,n\rightarrow \infty$ by Theorem~3.2. The second term also
approaches zero when $m,n\rightarrow \infty$ by Lemma~5.1(b) and
the condition (5.2). So proof is completed.\hfill $\Box$
\end{proof}

As a result we give the approximation order of
$L_{m,n}^{\alpha_{i},\beta_{j}},  j=1,2, m,n\in \mathbb{N},$ by
means of the weighted modulus of continuity.

\begin{theore}[\!]
For any $s>0$ and all $m,n\in N$ the inequality
\begin{equation}
\underset{\Vert f\Vert_{\rho }=1}{\sup }\left\{ \underset{\sqrt{
x^{2}+y^{2}}\leq s}{\sup }\vert L_{m,n}^{\alpha_{i},\beta
_{j}}(f;x,y) -f(x,y)\vert \right\} \leq c\underset{\Vert
f\Vert_{\rho }=1}{\sup }[ w_{\rho }(f,\delta ) ]
\end{equation}
holds for the linear positive operators
$\{L_{m,n}^{\alpha_{i},\beta_{j}}\}, j=1,2,$ defined on
$C_{\rho},$ where $\delta =\sqrt{L_{m,n}^{\alpha_{i},\beta_{j}}[
(t-x) ^{2}+(\tau -y) ^{2}]}$ and $c>0$ is a constant depending on
$s$.
\end{theore}

\begin{proof}
Since $L_{m,n}^{\alpha_{i},\beta_{j}}, j=1,2, m,n\in \mathbb{N},$
are linear positive operators, we have
\begin{align*}
&\vert L_{m,n}^{\alpha_{i},\beta_{j}}(f;x,y)
-f(x,y)\vert\\[.6pc]
&\quad\leq L_{m,n}^{\alpha_{i},\beta_{j}}(\vert f(t,\tau
)-f(x,y)\vert
;x,y)+\vert f(x,y)\vert (L_{m,n}^{\alpha_{i},\beta_{j}}(1;x,y) -1)\\[.6pc]
&\quad\leq L_{m,n}^{\alpha_{i},\beta_{j}}\left(\rho (x,y)w_{\rho
}\left[ f; \frac{\sqrt{(t-x) ^{2}+(\tau -y) ^{2}}}{\delta } \delta
\right] ;x,y\right),
\end{align*}
by Lemma~3.1. Using (2.5) we get
\begin{align}
&\vert L_{m,n}^{\alpha_{i},\beta_{j}}(f;x,y) -f(x,y)\vert\nonumber\\[.45pc]
&\quad\leq \rho (x,y)w_{\rho }(f;\delta
)L_{m,n}^{\alpha_{i},\beta_{j}}\left(1+\left[ \frac{\sqrt{(t-x)
^{2}+(\tau -y) ^{2}}}{\delta }\right]
;x,y\right)\nonumber\\[.6pc]
&\quad\leq \rho (x,y)w_{\rho }(f;\delta )L_{m,n}^{\alpha_{i},\beta
_{j}}\left(1+\left[ \frac{(t-x) ^{2}+(\tau -y) ^{2}}{\delta ^{2}}\right]; x,y\right)\nonumber\\[.6pc]
&\quad\leq \rho (x,y)w_{\rho }(f;\delta )L_{m,n}^{\alpha_{i},\beta
_{j}}(\rho ;x,y) \!+\!\frac{1}{\delta
^{2}}L_{m,n}^{\alpha_{i},\beta _{j}}([ (t\!-\!x) ^{2}+(\tau\!-\!y)
^{2}] ;x,y).
\end{align}
Since $(t-x) ^{2}+(\tau -y) ^{2}\in C_{\rho },$ (5.5) gives
\begin{align*}
\vert L_{m,n}^{\alpha_{i},\beta_{j}}(f;x,y) -f(x,y)\vert &\leq
\underset{\sqrt{x^{2}+y^{2}}\leq \ s}{\sup } c^{2}w_{\rho
}(f;\delta )\Vert L_{m,n}^{\alpha_{i},\beta_{j}}(\rho ;x,y)
\Vert_{\rho}\\[.6pc]
&\quad\,+\frac{1}{\delta ^{2}}\Vert L_{m,n}^{\alpha_{i},\beta
_{j}}((t-x) ^{2}+(\tau -y) ^{2};x,y) \Vert_{\rho },
\end{align*}
where $c=\sup_{\sqrt{x^{2}+y^{2}}\leq \ s}\rho (x,y).$ $\Vert
L_{m,n}^{\alpha_{i},\beta_{j}}(\rho ;x,y) \Vert_{\rho }$ is
bounded since
\begin{equation*}
\Vert L_{m,n}^{\alpha_{i},\beta_{j}}(\rho ;x,y) \Vert_{\rho
}=\Vert L_{m,n}^{\alpha_{i},\beta_{j}}(\rho ;x,y) \Vert_{C_{\rho
}\rightarrow B_{\rho }}
\end{equation*}
which is uniformly bounded, for all $m,n\in \mathbb{N},$ by
Lemma~5.1. From (5.5) and (5.6) we arrive at
\begin{equation*}
\underset{\sqrt{x^{2}+y^{2}}\leq \ s}{\sup }\vert L_{m,n}^{\alpha
_{i},\beta_{j}}(f;x,y) -f(x,y)\vert \leq c^{2}(1+M) w_{\rho
}(f;\delta ),
\end{equation*}%
which implies that
\begin{align*}
&\underset{\Vert f\Vert_{\rho }=1}{\sup }\left\{ \underset{\sqrt{
x^{2}+y^{2}}\leq s}{\sup }\vert L_{m,n}^{\alpha_{i},\beta
_{j}}(f;x,y) -f(x,y)\vert \right\}\\[.65pc]
&\quad\leq c^{2}(1+M) \underset{\Vert f\Vert_{\rho }=1}{\sup }[
w_{\rho }(f,\delta )],
\end{align*}
where $\delta =\sqrt{\Vert L_{m,n}^{\alpha_{i},\beta_{j}}([ (t-x)
^{2}+(\tau -y) ^{2}] ;x,y) \Vert_{\rho }}.$ Last inequality gives
(5.4), which completes the proof.\hfill $\Box$
\end{proof}


\begin{thebibliography}{99}
\bibitem{1} Alt\i n~A, Do\u{g}ru~O and Ta\c{s}delen~F, The generalization of
Meyer--K\"{o}nig and Zeller operators by generating functions,
{\it J. Math. Anal. Appl.} (in print)

\bibitem{2} Baskakov~V~A, On a construction of converging sequences of
linear positive operators, {\it Studies of Modern Problems of
Constructive Theory of Functions} (1961) 314--318

\bibitem{3} Becker~M, Global approximation theorems for Szazs--Mirakyan and
Baskakov operators in polynomial weight spaces, {\it Indiana Univ.
Math. J.} {\bf 27(1)} (1978) 127--142

\bibitem{4} Bernstein~S~N, D\'{e}monstration du th\'{e}or\'{e}me de
Weierstrass fond\'{e}e sur la calcules des probabilites, {\it
Comm. Soc. Math. Charkow S\'{e}r.} {\bf 13(2)} (1912) 1--2

\bibitem{5} B\"{u}y\"{u}kyaz\i c\i~$\dot{\rm I}$ and $\dot{\rm I}$bikli~E, The approximation
properties of generalized Bernstein polynomials of two variables,
{\it Appl. Math. Comput.} {\bf 156(2)} (2004) 367--380

\bibitem{6} Callahan~J, Advanced Calculus, Lecture Notes (USA: Smith College)

\bibitem{7} Zhou~D~X, Weighted approximation by Sz\'{a}sz--Mirakjan
operators, {\it J. Approx. Theory} {\bf 76(3)} (1994) 393--402

\bibitem{8} Do\u{g}ru~O, Weighted approximation of continuous functions on
the all positive axis by modified linear positive operators, {\it
Int. J. Comput. Numerical Anal. Appl.} {\bf 1(2)} (2002) 135--147

\bibitem{9} Do\u{g}ru~O, \"{O}zarslan~M~A and Ta\c{s}delen~F, On positive
operators involving a certain class of generating functions, {\it
Studia Scientiarum Mathematicarum Hungarica} {\bf 41(4)} (2004)
415--429

\bibitem{10} Gadjiev~A~D, Weighted approximation of continuous functions
by positive linear operators on the whole real axis, (Russian)
{\it Izv. Akad. Nauk Azerba\u\i d\v zan. SSR Ser. Fiz.-Tehn. Mat.
Nauk} {\bf 5} (1975) 41--45

\bibitem{11} Gad\v ziev~A~D, Positive linear
operators in weighted spaces of functions of several variables,
{\it Izv. Akad. Nauk Azerba\u\i dzhan. SSR Ser. Fiz.-Tekhn. Mat.
Nauk 1} {\bf 4} (1980) 32--37

\bibitem{12} Hermann~T, On the Sz\'{a}sz-Mirakian operator, {\it Acta Math. Acad.
Sci. Hungar.} {\bf 32(1--2)} (1978) 163--173

\bibitem{13} Ispir~N and Atakut~\c{C}, Approximation by modified
Szasz-Mirakjan operators on weighted spaces, {\it Proc. Indian
Acad. Sci. (Math. Sci.)} {\bf 112(4)} (2002) 571--578

\bibitem{14} Kirov~G~H, A generalization of the Bernstein
polynomials, {\it Math. Balkanica (N.S.)} {\bf 6(2)} (1992)
147--153

\bibitem{15} Kirov~G H and Popova~L, A generalization of the linear positive
operators, {\it Math. Balkanica (N.S.)} {\bf 7(2)} (1993) 149--162

\bibitem{16} Stancu~D~D, Approximation of functions by a new class of
linear polynomial operators, {\it Rev. Roumaine Math. Pures Appl.}
{\bf 13} (1968) 1173--1194

\bibitem{17} Stancu~D~D, A new class of uniform approximating polynomial
operators in two and several variables, {\it Proceedings of the
Conference on the Constructive Theory of Functions (Approximation
Theory)} (Budapest, 1969) pp.~443--455 (Budapest: Akad\'{e}miai
Kiad\'{o}) (1972)

\bibitem{18} Walczak~Z, Approximation of functions of two variables by some
linear positive operators, {\it Acta. Math. Univ. Comenianae} {\bf
LXXIV(1)} (2005) 37--48

\bibitem{19} Volkov~V~I, On the convergence of sequences of linear
positive operators in the space of continuous functions of two
variables, (Russian) {\it Dokl. Akad. Nauk SSSR (N.S.)} {\bf 115}
(1957) 17--19
\end{thebibliography}
\end{document}